\documentclass[times,polish,12pt]{article}
\usepackage[utf8]{inputenc}
\usepackage[T1]{fontenc}
\usepackage{amsmath,amssymb,times}
\usepackage{graphicx,colordvi,color}

\usepackage{url}

\begin{document}

\title{Inverted solutions of KdV-type and Gardner equations}

\author{Anna Karczewska\\ \small
Institute of Mathematics, Faculty of Mathematics, Computer Science and Econometrics\\ \small University of Zielona G\'ora, Szafrana 4a, 65-246 Zielona G\'ora, Poland\\ \small
A.Karczewska@wmie.uz.zgora.pl\\
\large Piotr Rozmej\\ \small
Institute of Physics, Faculty of Physics and Astronomy \\ \small
University of Zielona G\'ora, Szafrana 4a, 65-246 Zielona G\'ora, Poland\\ \small
P.Rozmej@if.uz.zgora.pl}

\date{\today} 
\maketitle

\begin{abstract}
In most of the studies concerning nonlinear wave equations of Korteweg-de Vries type, the authors focus on waves of elevation. Such waves have general form ~$u_{\text{u}}(x,t)=A f(x-vt)$, where ~$A>0$. In this communication we show that if ~$u_{\text{up}}(x,t)=A f(x-vt)$ is the solution of a given nonlinear equation, then $u_{\text{down}}(x,t)=-A f(x-vt)$, that is, an inverted wave is the solution of the same equation, but with changed sign of the parameter ~$\alpha$. This property is common for KdV, extended KdV, fifth-order KdV, Gardner equations, and generalizations for cases with an uneven bottom.
\end{abstract}


\section{Flat bottom case}\label{s1}
Begin with the simplest cases, when the bottom of the basin is flat. Because we will later also consider the cases with an uneven bottom, we will use equations corresponding to the fixed reference frame. We will also show that all properties discussed by us in this communication that are valid in a fixed reference frame remain valid in a frame moving with a constant velocity.

\subsection{Korteweg-de Vries equation}

 Korteweg-de Vries equation has in the fixed frame the following form (subscripts stand for partial derivatives)
\begin{equation} \label{kdv}	
u_t+u_x +\frac{3}{2}\alpha u u_x+\frac{1}{6}\beta u_{xxx} =0.
\end{equation}
Recall that the KdV equation applies for long surface waves with small amplitude, because it was derived under assumptions that two parameters 
\begin{equation} \label{lbeDef}
\alpha =\frac{A}{H}\qquad \mbox{and} \qquad \beta= \left(\frac{H}{L}\right)^{2}
\end{equation} 
are small and of the same order ($A$ denotes the wave amplitude, $H$ - the basin depth, and $L$ - the wavelength).

KdV equation (\ref{kdv}) has analytic $n$-soliton solutions and two kinds of periodic solutions, cnoidal and superposition ones. Single soliton solutions and periodic ones have the form of ``travelling waves'', that is, the waves with a fixed profile, moving with a constant speed. Such waves are described by the following functions
\begin{equation}\label{trwav}
u(x,t) = A\, f[B(x-vt)]+D.
\end{equation}
For single soliton solution one has (denoting $\xi=x-vt$)
\begin{equation}\label{1scoef}
f(\xi)= \text{sech}(B\xi)^{2}, \qquad B=\sqrt{\frac{3\alpha}{4\beta}A}, \qquad v=1+\frac{\alpha}{2}A, \qquad D=0,
\end{equation}
where $\text{sech}(\xi)\!=\!1/\text{cosh}(\xi)$ denotes the hyperbolic secant.
For cnoidal solutions to (\ref{kdv}) one has~(see~\cite{IKRR})
\begin{align}\label{percoef}
f(\xi)= \text{cn}^{2}(\xi,m),& \qquad B=\sqrt{\frac{3\alpha}{4\beta}\frac{A}{m}}, \qquad v=1+\frac{\alpha}{2}\frac{A}{m} \left(2-m-3\frac{E(m)}{K(m)} \right), \\ \nonumber \mbox{and} \quad & \qquad D=- \frac{A}{m} \left(\frac{E(m)}{K(m)} +m-1\right).
\end{align}
In (\ref{percoef}), $\text{cn}(\xi)$ is the Jacobi elliptic function, $E(m)$ is the complete elliptic integral,  $K(m)$  is the complete elliptic integral of the first kind, and $m$ is so-called elliptic parameter. The constant $D$ provides volume conservation (the volume of the elevated part of the fluid is equal to the volume of its lowered part).

Only since 2013 have superposition solutions in the following form been known (see  \cite{Khare})
\begin{equation} \label{Rsup}
u(x,t) = A\, f(\xi)+D=\frac{A}{2} \left[ \text{dn}^{2}(B\xi,m)\pm\sqrt{m}\,\text{cn}(B\xi,m)\text{dn}(B\xi,m)\right]+D, 
\end{equation}
where $\text{dn}(B\xi,m)$ is also one of Jacobi elliptic functions. The profiles of superposition solutions (\ref{Rsup}) are similar to profiles of cnoidal ones for the same $m$, but the velocity $v=1+\frac{\alpha A}{8}\left[5-m-6\frac{E(m)}{K(m)} \right]$ is slightly different \cite{RK,RKI,KRbook}.

The KdV equation is often considered in a moving frame. For instance, applying a transformation $\hat{x}=x-t, ~\hat{t}=t$ to (\ref{kdv}) one obtains the equation (where signs ˆ are already omitted)
\begin{equation} \label{kdvr} u_t+\frac{3}{2}\alpha u u_x+\frac{1}{6}\beta u_{xxx} =0. \end{equation} 
The form of solutions (\ref{trwav}) to the equation (\ref{kdvr}) and the coefficients given by (\ref{1scoef})-(\ref{percoef}) remain the same, only $\hat{v}=v-1$. So, the transformation to a moving frame $\hat{x}=x-t, ~\hat{t}=t$ removes the term $u_{x}$ from~(\ref{kdv}). The same term is removed when such transformation is applied to the extended KdV, fifth-order KdV and Gardner equations.

The literature usually considers solutions to the equation (\ref{kdv}) with $A>0$. Do functions with $A<0$ also satisfy the KdV equation? Suppose that the function (\ref{trwav}) satisfies the equation (\ref{kdv}). Therefore, the following occurs
\begin{equation} \label{kdvOK}
A f_{t} + A f_{x}+ \frac{3}{2} \alpha A^{2} f\,f_{x}+\frac{1}{6}\beta f_{xxx} =0.
\end{equation}
Inversion $u\to -u$ entails  $A \to -A$ and, by definition of $\alpha$, $\alpha \to -\alpha$, but leaves $B$ and $v$ unchanged. 
Thus, if $u$ given by (\ref{trwav}) with $A>0$ satisfies the equation (\ref{kdv}) with $\alpha>0$, then $u'=-u$ satisfies equation (\ref{kdv}) with $\alpha<0$ because substituting $u'=-A f$ into $u_t+u_x +\frac{3}{2}(-\alpha) u u_x+\frac{1}{6}\beta u_{xxx} $ merely changes the sign of the equation (\ref{kdvOK})
\begin{equation} \label{kdvOK-}
-A f_{t} - A f_{x}- \frac{3}{2} \alpha A^{2} f\,f_{x}-\frac{1}{6}\beta f_{xxx} 
= - (A f_{t} + A f_{x}+ \frac{3}{2} \alpha A^{2} f\,f_{x}+\frac{1}{6}\beta f_{xxx}) =0.
\end{equation}
The same property holds for periodic solutions, both cnoidal and superposition ones,  (\ref{percoef}) and (\ref{Rsup}), respectively.

The $n$-soliton solutions are given in terms of nonlinear superpositions, see e.g., \cite{Brauer}, where solutions with $n=2,3$ are presented in moving frame. For $n=2$ still, it is relatively easy to trace the consequences of changing the sign of the solution. Let us denote by $A_{2}>A_{1}>0$ the amplitudes of the higher and lower solitons. 
In the fixed frame, we introduce the following notations
\begin{equation} \label{Barg2}
\Theta_i(x,t)=\sqrt{\frac{3\,\alpha}{4\,\beta} A_i} 
   \left[x- t \left(1+\frac{\alpha}{2}  A_i\right)\right].
\end{equation}
Then, $2$-soliton solutions to (\ref{kdv}) have the following form \cite{KRbook}
\begin{align}\label{WarPoczZmBez2Sol}u(x,t) & =\frac{(A_2-A_1) \left(A_1 \,\text{sech}^2\left[\Theta_1(x,t)\right]+A_2\,\text{csch}^2\left[\Theta_2(x,t)\right]\right)}{\left(\sqrt{A_1} \tanh \left[\Theta_1(x,t)\right]-\sqrt{A_2}  \coth \left[\Theta_2(x,t)\right]\right)^2} .
\end{align}
We know that the function $u(x,t)$ given by (\ref{WarPoczZmBez2Sol}), for $A_{2}>A_{1}>0$ satisfies the KdV equation (\ref{kdv}) in which $\alpha>0$. 
Let us check how the solution (\ref{WarPoczZmBez2Sol}) will look like if we change the signs of the amplitudes of the two solitons and the sign of $\alpha$, i.e., substitute $A_{i}\to -A_{i}$ and $\alpha\to -\alpha$. Since these changes do not alter the arguments of $\Theta_{i}$, we obtain
\begin{align}\label{WarPoczZm}
u'(x,t) & =\frac{(-A_2+A_1) \left(-A_1 \,\text{sech}^2\left[\Theta_1(x,t)\right]-A_2\,\text{csch}^2\left[\Theta_2(x,t)\right]\right)}{\left(\sqrt{-A_1} \tanh \left[\Theta_1(x,t)\right]-\sqrt{-A_2}  \coth \left[\Theta_2(x,t)\right]\right)^2} 
\\ & = \frac{(A_2-A_1) \left(A_1 \,\text{sech}^2\left[\Theta_1(x,t)\right]+A_2\,\text{csch}^2\left[\Theta_2(x,t)\right]\right)}{(i)^{2}\left(\sqrt{A_1} \tanh \left[\Theta_1(x,t)\right]-\sqrt{A_2}  \coth \left[\Theta_2(x,t)\right]\right)^2} =-u(x,t). \nonumber
\end{align}
Indeed, the function (\ref{WarPoczZm}) represents a 2-soliton solution inverted with respect to (\ref{WarPoczZmBez2Sol}). It satisfies the equation (\ref{kdv}) in which $\alpha\to-\alpha<0$, which can be verified by direct calculus, since
$$ 0= -\left(u_t+u_x +\frac{3}{2}\alpha u u_x+\frac{1}{6}\beta u_{xxx}  \right) =
u'_t+u'_x +\frac{3}{2}(-\alpha) u' u'_x+\frac{1}{6}\beta u'_{xxx}. $$

The explicit form of $3$-soliton solution is more complicated \cite{KRbook,Brauer}. Assume $A_{3}>A_{2}>A_{1}>0$. Denote
\begin{align}\label{3SolA}
X_1(x,t) & := -\frac{2(A_1-A_2)\left(A_1\,\text{sech}^2\left[\Theta_1(x,t)\right]+A_2\, \text{csch}^2\left[\Theta_2(x,t)\right]\right)}{\left(\sqrt{2 A_1} \tanh \left[\Theta_1(x,t)\right] - \sqrt{2 A_2}  \coth \left[\Theta_2(x,t)\right] \right)^2},
\\ \label{3SolB}
X_2(x,t) & := \frac{2(-A_1+A_3)\left(-A_1\,\text{sech}^2\left[\Theta_1(x,t)\right]+A_3\, \text{sech}^2\left[\Theta_3(x,t)\right]\right)}{\left(-\sqrt{2 A_1} \tanh \left[\Theta_1(x,t)\right] + \sqrt{2 A_3}  \tanh \left[\Theta_3(x,t)\right] \right)^2},
\\ \label{3SolC}
X_3(x,t) & := \frac{2(A_1-A_2)}{\sqrt{2A_1}\tanh \left[\Theta_1(x,t)\right]
-\sqrt{2A_2}\coth \left[\Theta_2(x,t)\right]},
\\ \label{3SolD}
X_4(x,t) & := \frac{2(-A_1+A_3)}{-\sqrt{2A_1}\tanh \left[\Theta_1(x,t)\right]+
\sqrt{2A_3}\tanh \left[\Theta_3(x,t)\right]}.
\end{align}
Then 3-soliton solution is expressed with  (\ref{3SolA})-(\ref{3SolD}) as
\index{KdV!3-soliton solution}
\begin{equation}\label{3Sol}
u(x,t)= A_1\, \text{sech}^2\left[\Theta_1(x,t)\right] - 2(A_2-A_3) \frac 
{X_1(x,t)+X_2(x,t) }{\left(X_3(x,t)-X_4(x,t)\right)^2}.
\end{equation}
With the same arguments, we see that the inversion $A_{i}\to -A_{i}, ~i=1,2,3$ implies $u'\to -u$ for the formula (\ref{3Sol}). Therefore, the inverted 3-soliton solution fulfilles the KdV equation with $\alpha'=\alpha$.

Analogously, similar arguments hold for $n$-soliton solutions and for KdV equation in moving frame (\ref{kdvr}).

\subsection{Extended KdV equation (KdV2)}

The extended KdV equation (KdV2), first derived by Marchant and Smyth \cite{MS90}, has in 
the fixed reference frame the following form
\begin{align} \label{eta2}
u_t+u_{x} & + \frac{3}{2} \alpha u u_x +\frac{1}{6} \beta u_{xxx} 
-\frac{3}{8} \alpha^2 u^2 u_x 
+ \alpha\beta \left(\frac{23}{24}u_x u_{xx}+\frac{5}{12}u u_{xxx}\right)+\frac{19}{360} \beta^2 u_{xxxxx}=0. 
\end{align}

Both the KdV equation (\ref{kdv}) and the KdV2 equation (\ref{eta2}) are derived from the model of an ideal fluid (incompressible and inviscid) with an irrotational motion under the assumption that $\alpha=O(\beta)$, that is, the parameters $\alpha,\beta$ are small and of the same order. Then, applying the perturbation technique, one obtains the KdV equation when the approach is limited to first order in small parameters and the extended KdV equation (KdV2) when the derivations are extended to second order.

The equation (\ref{eta2}), unlike the KdV equation (\ref{kdv}), is nonintegrable. Despite this fact, it has analytic single-soliton solutions \cite{KRI}, as well as periodic
 cnoidal solutions \cite{IKRR}  and superposition solutions \cite{RK,RKI}.  However, multi-soliton solution to (\ref{eta2}) do not exist \cite{KRact}.
 
 The analytic single-soliton and periodic solutions of the KdV2 equation have the same functional form as the analogous solutions of the KdV equation but the corresponding coefficients $B,v$ are slightly different, and the coefficient $A$ is determined by the parameters of the equation, i.e., by $\alpha,\beta$. It is because the KdV2 equation imposes one more condition on the coefficients $A,B,v$, determining the solution, than the KdV equation \cite{IKRR,RK,RKI,KRbook,KRI}. For KdV, the set of these coefficients has one degree of freedom, which, for fixed $\alpha,\beta$, allows solitons of different heights to exist and hence admits multi-soliton solutions. For KdV2 with fixed $\alpha,\beta$, there exists a soliton with the only possible height, so multi-soliton solutions do not exist  \cite{KRact}.

If the function $u(x,t)$ satisfies the KdV2 equation (\ref{eta2}), it is easy to see that the inverted function $u'=-u$ satisfies the equation (\ref{eta2}) with the sign of $\alpha$ changed $\alpha'=-\alpha$~ because
\begin{align} \label{2s1}
& ~u'_t +u'_x +\frac{3}{2}\alpha' u' u'_x+\frac{1}{6}\beta u'_{xxx}-\frac{3}{8} \alpha'^2 u'^2 u'_x + \alpha'\beta \left(\frac{23}{24}u'_x u'_{xx}+\frac{5}{12}u' u'_{xxx}\right)+\frac{19}{360} \beta^2 u'_{xxxxx} \nonumber \\
= &-u_t - u_{x}-\frac{3}{2}\alpha u u_x-\frac{1}{6}\beta u_{xxx}+\frac{3}{8} \alpha^2 u^2 u_x - \alpha\beta \left(\frac{23}{24}u_x u_{xx}+\frac{5}{12}u u_{xxx}\right)-\frac{19}{360} \beta^2 u_{xxxxx} \\
= &-\left(u_t +u_{x}+ \frac{3}{2} \alpha u u_x +\frac{1}{6} \beta u_{xxx} 
-\frac{3}{8} \alpha^2 u^2 u_x 
+ \alpha\beta \left(\frac{23}{24}u_x u_{xx}+\frac{5}{12}u u_{xxx}\right)+\frac{19}{360} \beta^2 u_{xxxxx} \right)=0. \nonumber
\end{align}

This conclusion remains valid also in a moving frame of reference.

\subsection{Fifth-order KdV equation}

The equations KdV and KdV2 were derived with complete neglect of the pressure exerted by the surface tension of the deformed liquid surface. The terms derived from the surface tension include a coefficient called the Bond number~ $\tau=\frac{T}{\varrho g H^{2}}>0$, where $T$~ is the surface tension coefficient of the liquid, $\varrho$~ is its density, and $H$~ is the depth of the reservoir. For water, at a depth of $H\geq 1$m, we have $T \le 10^{-7}$, so the effect of surface tension can be completely neglected. However, for thin liquid layers, when $H$~ is on the order of millimeters, the terms derived from the surface tension become significant \cite{KRcnsns}. 

Taking into account the terms originating from surface tension in the Euler equations, assuming $\alpha=O(\beta^{2})$, and applying perturbation approach up to second order in small parameters, one derives so-called \emph{fifth-order KdV]} equation (see \cite{HS88,BS13})
\begin{align}\label{5kdvQ}
u_t+u_x  & + \frac{3}{2}\alpha u u_x +\beta \frac{1-3\tau}{6} u_{xxx}  +\beta^2 \frac{19-30\tau-45\tau^2}{360}u_{xxxxx}   
=0.\end{align}
It is well known, see, e.g.~\cite{Dey96,Bri02}, that the fifth-order KdV equation has a  soliton solution in the form
\begin{equation} \label{5sol1}
u(x,t) = A\, \text{sech}^4[B(x-vt)].
\end{equation}

Here it is also easy to see that if $u$ satisfies the equation (\ref{5kdvQ}) then the function $u'=-u$ satisfies the equation (\ref{5kdvQ}) with $\alpha'=-\alpha$.

\subsection{Gardner equation}
If one assumes that $\beta=O(\alpha^{2})$ then the so-called Gardner equation \cite{BS13} is derived from Euler's equations with second order accuracy in small parameters
\begin{align} \label{Gard}
u_t+u_x +\frac{3}{2} \alpha  u u_x & - \frac{3}{8} \alpha^2 u^2 u_x + \frac{1-3\tau}{6} \beta\, u_{xxx}  =0 .
\end{align}
It is a well-known fact that for the Gardner equation (\ref{Gard}) there exists one-parameter family of analytic solutions in the form \cite{GPT99,OPSS15}
\begin{equation} \label{1sOPSS}
u(x,t)= \frac{A}{1+ B\, \textrm{cosh}[(x-v\,t)/\Delta]}.
\end{equation}
The functions (\ref{1sOPSS}) also represent solitons with fixed profiles, moving with a constant speed. Their shapes can be very different, from bell-shaped solitons to table-top ones, depending on values of equation parameters \cite{GPT99,OPSS15}.

The equation (\ref{Gard}) imposes three conditions on coefficients $A,B,v,\Delta$ of solutions. So, three of them can be expressed as functions of the single one. Choosing $\Delta$ as the independent parameter, one obtains the following relations (denoting $\frac{(1-3\tau)}{6}\beta=\beta'$)
\begin{equation} \label{Gpar}
A=\frac{4\,\beta'}{\alpha\,\Delta^2}, \qquad B= \pm\sqrt{1-\frac{\beta'}{\Delta^2}}, \qquad v=1+\frac{\beta'}{\Delta^2}.
\end{equation}
In (\ref{Gpar}), $B,v$ do not depend on $\alpha$. Then inverted function $u'=-u$ requires 
$\alpha'=-\alpha$ with $B,v,\Delta$ unchanged. 
In this case it is also seen that if $u$ satisfies equation (22) then the function $u'=-u$ satisfies equation (22), where $\alpha' = -\alpha$.

{\bf Conclusion 1:  The inverted solutions of the KdV, KdV2, fifth-order KdV, and Gardner equations satisfy these equations with a negative $\alpha$ parameter.}

\section{Uneven bottom}

Back in 2014, in the work \cite{KRI}, we started a systematic study on the extension of the KdV and KdV2 equations to the more general case when the bottom of the tank is described by a continuous bounded function $h(x)$. To this end, we introduced a new small parameter $\delta= \frac{a_{h}}{H}$, where $a_{h}$ denotes the amplitude of the bottom profile.

However, it is only in the work of \cite{KRcnsns} that we have fully realized this goal. Moreover, we have shown that the new term in these nonlinear wave equations is universal, i.e., it is the same for all four equations: KdV, KdV2, fifth-order KdV, and Gardner. The same constraint holds in all these cases. Namely, the system of Boussinesq equations can be reduced to a single wave equation only if the bottom function satisfies the condition $h_{xx}=0$, i.e., $h(x)$ is a piecewise linear function. In such cases, the generalized KdV, KdV2, fifth-order KdV and Gardner equations take the form
\begin{equation} \label{gen}
EQ -\frac{1}{4} \delta (2 h u_x+h_x u) =0,
\end{equation}
where $EQ$ stands for left hand side of equations (\ref{kdv}), (\ref{eta2}), (\ref{5kdvQ}) or (\ref{Gard}). 

The equations (\ref{gen}) are nonintegrable. Based on the analysis in section \ref{s1}, we know that if $u$ satisfies the equation $EQ=0$, then $u'=-u$ satisfies the equation $EQ'=0$ in which $\alpha'= -\alpha$, because then every term of the equation will change sign. It can be seen that changing $u\to -u$ will also change the sign of the entire term $-\frac{1}{4} \delta (2 h u_x+h_x u)$, if $\delta$ remains unchanged. Thus, if the function $u(x,t)$ is a solution to one of the generalized equations (\ref{gen}) (i.e., the KdV, KdV2, 5th-order KdV, or Gardner equations), then the inverted function 
$-u(x,t)$ is the solution of the corresponding one of these equations, in which $\alpha$ is replaced by $-\alpha$.

{\bf Conclusion 2:  The inverted solutions of the KdV, KdV2, fifth-order KdV, and Gardner equations generalized for uneven bottom satisfy these equations with a negative $\alpha$ parameter.}


\end{document}